\documentclass[12pt]{article}
\usepackage{amsmath,amssymb,amsthm}
\usepackage{graphicx,authblk}
\usepackage{subfigure}
\usepackage{amsmath,amsfonts,amssymb,amscd}
\usepackage{indentfirst,graphicx,epstopdf}
\usepackage{graphicx,psfrag,ifpdf,enumerate}
\usepackage{caption}
\usepackage{amsthm, amsfonts, color}
\usepackage[bookmarksnumbered, plainpages]{hyperref}
\usepackage[numbers,sort&compress]{natbib}
\usepackage[noend]{algpseudocode}
\usepackage{algorithmicx,algorithm}
\usepackage{psfrag}
\usepackage{color}
\usepackage{enumerate}
\usepackage{epstopdf}
\input{epsf}

\newtheorem{theorem}{Theorem}[section]
\newtheorem{conj}[theorem]{Conjecture}

\theoremstyle{definition}

\title{\Large \bf Short rainbow cycles in edge-colored graphs}

\author[a]{Xiaozheng Chen}
\author[a]{Shanshan Guo}
\author[a]{Fei Huang\thanks{Corresponding author: Fei Huang. Email: hf@zzu.edu.cn}}
\affil[a]{School of Mathematics and Statistics, Zhengzhou University, Zhengzhou, Henan, China \authorcr \it{ Email: cxz@zzu.edu.cn; 15738385820@163.com; hf@zzu.edu.cn.}}

\date{}
\begin{document}
\maketitle

\begin{abstract}
A famous conjecture of Caccetta and H\"{a}ggkvist (CHC) states that
a directed graph $D$ with $n$ vertices and minimum outdegree at least $r$
has a directed cycle of length at most $\lceil \frac{n}{r}\rceil$.
In 2017, Aharoni proposed the following generalization:
an edge-colored graph $G$ with $n$ vertices, $n$ color classes of size at least $r$
has a rainbow cycle of length at most $\lceil \frac{n}{r}\rceil$.
Since CHC can be seen as the case of Aharoni's Conjecture:
color classes in the color partition are monochromatic stars centered at distinct vertices,
one way to study Aharoni's Conjecture is to structure the color classes
as each color class is either a star, a triangle or contains a matching of size 2.
Guo improved the upper bound in Aharoni's Conjecture to $O(\log n)$ in some mixed cases when the color classes are not necessarily stars.
In this paper, we extend Guo's results. Our main result is as follows: Let $G$  an edge-colored graph on $n$ vertices and $n$ color classes,
if at least $\alpha n$ color classes  are either
a matching of size 2 or a triangle for $\alpha >\frac{1}{2}$,
then $G$ contains a rainbow cycle of length $O(\log n)$. We also prove that the $\log n$ bound is the right order of magnitude.
\end{abstract}
\section{Introduction}

The \emph{directed girth} of a directed graph $D$ denoted by $dgirth(D)$ is the smallest length of a directed cycle in $D$.
The following famous conjecture on directed girth is proposed by  Caccetta and H\"{a}ggkvist \cite{CH}.

\begin{conj}[Caccetta-H\"{a}ggkvist Conjecture \cite{CH}]\label{conj:CH}
Let $D$ be a digraph on $n$ vertices.
For every positive integer $r$, if the minimum outdegree of $D$ is at least $r$,
then $dgirth(D)\leq \lceil\frac{n}{r}\rceil$.
\end{conj}

Despite the Caccetta-H\"{a}ggkvist Conjecture remains open,
it was proved that the conjecture is asymptotically true:
$dgirth(D)\leq \frac{n}{r}+73$ for all $D$ by Shen \cite{S}.

A graph $G=\left(V(G),E(G)\right)$ is edge-colored if there is a mapping $c:E(G)\rightarrow \mathbb{N}$.
Let $G$ be an edge-colored graph with a \emph{color partition} $\{F_1,\ldots,F_p\}$ of $E(G)$,
where $F_i=\{e\in E(G):c(e)=i\}$ for $1\leq i\leq p$.
The sets $F_i$ are called \emph{color classes}, the indices $i\in \{1,\ldots,p\}$ are called \emph{colors}.
A subgraph $H$ of $G$ is called \emph{rainbow} if distinct edges of $H$ come from distinct color classes.
The \emph{rainbow girth} of an edge-colored graph $G$ denoted by $rgirth(G)$ is the smallest length of a rainbow cycle in $G$.
In 2017, Aharoni proposed a conjecture which can be seen as a generalization of Caccetta-H\"{a}ggkvist Conjecture.

\begin{conj}[Aharoni \cite{ADH}]\label{conj:A}
Let $G$ be an edge-colored graph on $n$ vertices.
For every positive integer $r$, if $G$ has $n$ color classes and each color class is of size at least $r$,
then $rgirth(G)\leq \lceil\frac{n}{r}\rceil$.
\end{conj}

There are also some asymptotical results on Conjecture \ref{conj:A}:
(1) the conjecture is true if each color class has size at least $301r\log r$,
which was proved by Hompe et. al. \cite{HPPS};
(2) the conjecture is true if each color class has size at least $10^{11}r$,
which was proved by Hompe and Spirkl \cite{HP1};
(3) there exists a constant $c_r$ such that $rgirth(G)\leq \lceil\frac{n}{r}\rceil+c_r$, which was proved by Hompe and Huynh \cite{HH}.

Let $D$ be a digraph on $n$ vertices with minimum outdegree at least $r$.
To see Conjecture \ref{conj:A} implies Caccetta-H\"{a}ggkvist Conjecture,
we can get an edge-colored graph $G$ with color partition
$\{A^+_D(v)\mid v\in V(D)\}$ by the undirected underlying graph of $D$,
where $A^+_D(v)$ denotes the set of out arcs of $v$ in $D$.
It is easy to see that the cycle in $G$ is rainbow if and only if it is a directed cycle in the original digraph $D$.
Thus, Caccetta-H\"{a}ggkvist Conjecture is equivalent to the case of Conjecture \ref{conj:A}:
color classes in the color partition are
monochromatic stars centered at distinct vertices.

One may wonder what happens when not all the color classes are stars.
If a color class is not a star,
then it must be a triangle or
contain a matching of size 2.
In the extreme case when each color class is a matching of size 2,
Aharoni and Guo \cite{AG} proved the following result (in this paper, the logarithms are to base 2).

\begin{theorem}[\cite{AG}]
For any constant $\alpha>\frac{3\sqrt{6}}{8}$,
there exist an integer $N_0$ and an absolute constant $C$ such that
if $G$ is a simple edge-colored graph with $n$ vertices ($n\geq N_0$), $\alpha n$ colors , and each color class is a matching of size 2,
then $G$ contains a rainbow cycle of length at most $C\log n$.
\end{theorem}

 The other extreme case that each color class is
a triangle is also proved by Aharoni et. al. \cite{ABCGZ}.
They also showed that $\log n$ is the right order of magnitude and $\alpha>\frac{1}{2}$ is sharp.

\begin{theorem}[\cite{ABCGZ}]\label{thm:A}
For any constant $\alpha>\frac{1}{2}$,
there exists an absolute constant $C$ such that
if $G$ is a simple edge-colored graph with $n$ vertices, $\alpha n$ colors, and each color class is a triangle,
then $G$ contains a rainbow cycle of length at most $C\log n$.
\end{theorem}

Some mixed cases were investigated by Guo
and he showed that
$\alpha>\frac{1}{2}$ is sharp.

\begin{theorem}[\cite{G}]
There exist an integer $N_0$ and an absolute constant $C$ such that
if $G$ is a simple edge-colored graph with $n$ vertices ($n\geq N_0$), $n$ colors, and each color class is a matching of size 2 or a triangle,
then $G$ contains a rainbow cycle of length at most $C\log n$.
\end{theorem}

\begin{theorem}[\cite{G}]\label{thm:Guo2}
For any $\alpha>\frac{1}{2}$,
there exist an integer $N_0$ and an absolute constant $C$ such that
if $G$ is a simple edge-colored graph with $n$ vertices ($n\geq N_0$), $n$ colors,
and  at least $\alpha n$ color classes  are matchings of size 2,
then $G$ contains a rainbow cycle of length at most $C\log n$.
\end{theorem}

From the above results, we can see that if the structures of the color classes are not necessarily stars, the upper bound of the rainbow girth may be much smaller than the bound (linear of $n$) in  Conjecture \ref{conj:A}. On the other hand, we can construct an edge-colored graph such that the all the color classes are stars, and the rainbow girth reaches the upper bound in Conjecture \ref{conj:A}. Let $G$ be the edge-colored graph
with vertex set $V(G)=\{v_1,v_2,\ldots,v_{kr+1}\}$
such that $c(v_iv_{i+j})=i$ for $1\leq j\leq k$ (the subscript is taken module $kr+1$).
Apparently, $v_1v_{1+k}v_{1+2k}\ldots v_{1+rk}v_1$ is a shortest rainbow cycle of $G$, and so, the rainbow girth of $G$ is $r+1 =\lceil\frac{|V(G_1)|}{k}\rceil$. Thus, the proportion of the color classes of stars may effect the upper bound of rainbow girth.

In this paper, we prove the following result
which states that if there are more than $\frac{n}{2}$ color classes consisting of a triangle, or containing a matching of size 2,
then the rainbow girth of $G$ is at most $C\log n$.

\begin{theorem}\label{thm:main}
For any $\alpha>\frac{1}{2}$,
there exist an integer $N_0$ and an absolute constant $C$ such that
if $G$ is a simple edge-colored graph with $n$ vertices ($n\geq N_0$), $n$ colors,
and at least $\alpha n$ color classes  are either
a matching of size 2 or a triangle,
then $G$ contains a rainbow cycle of length at most $C\log n$.
\end{theorem}

Compared with Theorem \ref{thm:Guo2},
if an edge-colored graph contains three kinds of color classes:
matching of size 2, triangle, and star,
Theorem \ref{thm:main} can solve the case with more color classes of stars.
We can also give an example to imply
that $\alpha>\frac{1}{2}$ is sharp:
Assume $n$ is divisible by 6
and $H$ is an edge-colored graph
with $V(H)=\{v_{i,j}\mid i=1,\ldots,\frac{n}{6},j=1,\ldots,6\}$.
Let $\mathcal{F}_T\cup \mathcal{F}_M\cup \mathcal{F}_S$ be a color partition of $H$
where $\mathcal{F}_T=\{\{v_{i,1}v_{i,2},v_{i,2}v_{i,3},
v_{i,1}v_{i,3}\},\{v_{i,4}v_{i,5},v_{i,5}v_{i,6},
v_{i,4}v_{i,6}\}\mid i=1,\ldots,\frac{n}{6}\}$,
$\mathcal{F}_M=\{\{v_{i,2}v_{i,4},v_{i,3}v_{i,5}\}\mid i=1,\ldots,\frac{n}{6}\}$,
$\mathcal{F}_S=\{\{v_{i,2}v_{i,5}\},\{v_{i,3}v_{i,4}\},\{v_{i,6}v_{i+1,1}\}\mid i=1,\ldots,\frac{n}{6}\}$
(with the first subscript module $\frac{n}{6}$).
Then $|\mathcal{F}_T|+|\mathcal{F}_M|=|\mathcal{F}_S|$.
It can be verified that $H$ has a rainbow girth $\frac{2n}{3}$.

We  prove the following result which implies that the $\log n$ bound in Theorem \ref{thm:main} is
the right order of magnitude.

\begin{theorem}\label{thm:main3}
Let $G$ be an edge-colored graph with $n$ vertices.
Then there exists a constant $c>0$ such that
if each color class of $G$ is a matching of size 2,
or a triangle, or a single edge, then any rainbow cycle in $G$ has length at least $c \log n$.
\end{theorem}

This paper is organized as follows:
In section 2, we list some preliminaries to prove the main results.
In section 3, we give the proofs of Theorems \ref{thm:main} and \ref{thm:main3}.

\section{Preliminaries}
As in \cite{AG,G}, a key ingredient in the proof is a famous result
by Bollob\'{a}s and Szemer\'{e}di on the girth of sparse graphs.

\begin{theorem}[Bollob\'{a}s and Szemer\'{e}di \cite{BS}]\label{thm:BS}
For $n\geq 4$ and $k\geq 2$, every $n$-vertex graph with
$n+k$ edges has a girth at most
$$
\dfrac{2(n+k)}{3k}(\operatorname{log}k+\operatorname{loglog}k+4).$$
\end{theorem}

If an edge-colored graph $G$ has a subgraph $H$ of $G$
such that $|V(H)|\leq \beta n$ and $H$ has a rainbow edge set of size at least $(\beta+\gamma)n$,
then by Theorem \ref{thm:BS} there exists an absolute constant $C(\beta,\gamma)$
such that
$H$ contains a rainbow cycle of length at most $C(\beta,\gamma)\log n$.
The idea to prove Theorem \ref{thm:main}
is to find a subgraph $H$ satisfying the assumption above
by probabilistic methods.
We shall use the following well-known concentration inequalities.

\begin{theorem}[Markov's Inequality]\label{thm:Markov}
Let $X$ be a non-negative random variable.
For any $t>0$, we have $\mathbf{P}(X\geq t)\leq \frac{\mathbb{E}X}{t}$.
\end{theorem}

\begin{theorem}[Chernoff]\label{thm:Chernoff}
Let $X$ be a binomial random variable $Bin(n,p)$.
 For any $\sigma\geq0$, we have
$\mathbf{P}(X\geq (1+\sigma)\mathbb{E}X)\leq \exp \big(-\frac{\sigma^2\mathbb{E}X}{3}\big)$
and $\mathbf{P}(X\leq (1-\sigma)\mathbb{E}X)\leq \exp \big(-\frac{\sigma^2\mathbb{E}X}{2}\big)$.
\end{theorem}

\begin{theorem}[Chebyshev's Inequality]\label{thm:Chebyshev}
Let $X$ be a random variable. For any $\sigma>0$,
we have $\mathbf{P}(|X-\mathbb{E}X|>\sigma \mathbb{E}X )\leq \dfrac{VarX}{({\sigma\mathbb{E}X})^2}$.
\end{theorem}

\section{Main results}

\subsection{Proof of Theorem \ref{thm:main}}
To prove Theorem \ref{thm:main}, we give the following stronger result
which restricts the structures of color classes.

\begin{theorem}\label{thm:main2}
For any $\alpha>\frac{1}{2}$, there exist
an integer $N_0$ and an absolute constant $C$ such that the following hold.
Let  $G$ be a simple edge-colored graph with $n\geq N_0$ vertices and a color partition $\mathcal{F}=\left(F_1, \ldots, F_n\right)$.
If $\mathcal{F}=\mathcal{F}_M\cup \mathcal{F}_T\cup \mathcal{F}_S$, where

(1) every $F_{i} \in \mathcal{F}_M$ is a matching of size 2,

(2) every $F_{i} \in \mathcal{F}_T$ is a triangle,

(3) every $F_{i} \in \mathcal{F}_S$ is a single edge,

(4) $\left|\mathcal{F}_M\cup \mathcal{F}_T\right| \geq \alpha n$,\\
then $G$ contains a rainbow cycle of length at most $C\log n$.
\end{theorem}

Theorem \ref{thm:main2} implies Theorem \ref{thm:main} as
for a color class $F\in \mathcal{F}$, if it is a star of size at least 3
we can delete the edges in $F$ such that $F$ is a single edge,
if it is not a star
we can delete the edges in $F$ such that $F$ is either a triangle or a matching of size 2.

\begin{theorem}\label{7}
For any $\alpha>\frac{1}{2}$, there exist $\beta, \gamma>0$ and an integer $N_0$
such that for any $n\geq N_0$, given an $n$-vertex graph $G$ and an edge coloring of $G$ satisfying the assumption in Theorem \ref{thm:main2}, there exists a subset $H$ of $V(G)$ of size at most $\beta n$ such that $G[H]$ contains a rainbow edge set of size at least $(\beta+\gamma) n$.
\end{theorem}

\noindent{\bf Proof of Theorem \ref{7}.}
Since $\alpha>\frac{1}{2}$,
if $|\mathcal{F}_M|=o(n)$,
then there exists an integer $N_1$
such that for any $n\geq N_1$
we have $\alpha> \frac{1}{2}+\frac{|\mathcal{F}_M|}{n}$.
Then there exists an constant $\alpha'>\frac{1}{2}$
such that $|\mathcal{F}_T|\geq \alpha n-|\mathcal{F}_M|\geq \alpha'n$.
Hence the result follows by Theorem \ref{thm:A}
when $n\geq N_1$.
Similarly if $|\mathcal{F}_T|=o(n)$,
then there exists an integer $N_2$
such that
the result follows by Theorem \ref{thm:Guo2} when $n\geq N_2$.
Since $\alpha n\leq |\mathcal{F}_M|+|\mathcal{F}_T|\leq n$,
we may assume that  $|\mathcal{F}_M|=\Omega(n)$,
$|\mathcal{F}_T|=\Omega(n)$ and $n\geq \max\{N_1,N_2\}$.

Let $H$ be a random vertex subset of $V(G)$,
in which each vertex of $V(G)$ is included independently
  with probability $p$.
Let $X_i$ be the indicator random variable that
at least one edge in $F_i\in \mathcal{F}$ is contained in $G[H]$ and $
X_{\Theta}:=\sum_{F_i \in \mathcal{F}_{\Theta}} X_i
$, for $\Theta\in\{M,T,S\}$.
For any $\varepsilon>0$,
let $\mathfrak{A}_{\Theta}(\varepsilon)$
be the event that the number of color classes
in $\mathcal{F}_{\Theta}$ that have at least one edge contained in $G[H]$ is at most $(1-\varepsilon)\cdot\mathbb{E}X_{\Theta}$, for $\Theta\in\{M,T,S\}$.

Since each vertex is included in $H$ independently  with probability $p$,
we have
$$\mathbf{P}(X_i=1)=
\begin{cases}
2 p^2-p^4,&F_i\in \mathcal{F}_{M};\\
3p^2-2p^3,&F_i\in \mathcal{F}_{T};\\
p^2,&F_i\in \mathcal{F}_{S},\\
\end{cases}$$
and then
$$
\mathbb{E} X_{\Theta}=
\begin{cases}
|\mathcal{F}_{M}|\cdot(2 p^2-p^4),&{\Theta}=M;\\
|\mathcal{F}_{T}|\cdot(3p^2-2p^3),&{\Theta}=T;\\
|\mathcal{F}_{S}|\cdot p^2,&{\Theta}=S.\\
\end{cases}
$$
By Chebyshev's inequality,
for any $\varepsilon>0$,
we have
\begin{equation}\label{e:1}
\mathbf{P}(X_{\Theta}\leq (1-\varepsilon)\mathbb{E}X_{\Theta})\leq \dfrac{\operatorname{Var} X_{\Theta}}{{(\varepsilon\mathbb{E}X_{\Theta})}^2},
\text{\quad for $\Theta\in\{M,T,S\}$.}
\end{equation}
Since $
\operatorname{Var} X_{\Theta}=\mathbb{E} (X_{\Theta})^2-(\mathbb{E} X_{\Theta})^2=\sum_{F_i , F_j\in \mathcal{F}_{\Theta}}\left(\mathbb{E} (X_i X_j)-\mathbb{E} X_j \mathbb{E} X_j\right),
$
we can  estimate
the contribution of each pair $(F_i,F_j)$ for $F_i,F_j\in \mathcal{F}_\Theta$ to $\operatorname{Var} X_\Theta$.


\vspace{2mm}
\noindent{\bf Claim 1.} For any $\delta>0$,
we have
$0\leq \operatorname{Var} X_M\leq \delta(\mathbb{E} X_{M})^2+o((\mathbb{E} X_{M})^2)$.

\begin{proof}
For any $F_i,F_j\in \mathcal{F}_M$,
if the edges in $F_i$ and $F_j$ are vertex-disjoint,
then $X_i$ and $X_j$ are independent and $\mathbb{E}( X_i X_j)-\mathbb{E} X_i \mathbb{E} X_j=0$.
Thus, assume that $F_i$ and $F_j$ are joint.
Let $F_i=\{a, b\}$ and $F_j=\{c, d\}$.
Then there are totally three cases if $F_i$ and $F_j$ are joint.

(1) $F_i$ has only one edge that is joint with exactly one edge in $F_j$,
that is, $F_i\cup F_j$ has three components consisting of one 2-path and two independent edges.
W.l.o.g., let $a=x y, b=u v$ and $c=xz, d= st$.

(2) Each edge in $F_i$ is joint with exactly one edge in $F_j$,
that is, $F_i\cup F_j$ has two vertex-disjoint 2-paths.
W.l.o.g., let $a=x y, b=u v$ and $c=xz, d= us$.

(3) $F_i$ has at least one edge that is joint with
both edges in $F_j$.

\vspace{2mm}
For (1) and (2), we give the following table
in which $A$ denotes the event that the edges are in $G[H]$,
$B$ denotes the event that the edges are not in $G[H]$,
and
$P_k$ denotes the probability that both the events $A$ and $B$ happen in case ($k$), $k=1,2$.

$$\begin{tabular}{|p{3cm}|p{4cm}|p{4cm}|p{4cm}|}
\hline
$A$:&$B$:&$P_1(A\cap B)$&$P_2(A\cap B)$\\
belong to $G[H]$&not belong to $G[H]$&&\\
\hline
$a, c$&$b, d$&$p^3\left(1-p^2\right)^2$&$p^3\left((1-p)+p(1-p)^2\right)$\\ \hline
$a, d $&$b, c$&$p^4(1-p)\left(1-p^2\right)$&$p^4(1-p)^2$\\ \hline
$b, c$&$a, d$&$p^4(1-p)\left(1-p^2\right)$&$p^4(1-p)^2$\\ \hline
$b, d $&$a, c$&$p^4\left(p(1-p)^2+(1-p)\right)$&$p^3\left((1-p)+p(1-p)^2\right)$\\ \hline
$a, b, c$&$d$&$p^5\left(1-p^2\right)$&$p^5\left(1-p\right)$\\ \hline
$a, b, d$&$c$&$p^6(1-p)$&$p^5\left(1-p\right)$\\ \hline
$c, d, a$&$b$&$p^5\left(1-p^2\right)$&$p^5\left(1-p\right)$\\ \hline
$c, d, b$&$a$&$p^6(1-p)$&$p^5\left(1-p\right)$\\ \hline
$a, b, c, d$&-&$p^7$&$p^6$\\ \hline
\end{tabular} $$

Therefore, for case (1),
we have $$
\begin{aligned}
\mathbb{E} (X_i X_j)=&p^3\left(1-p^2\right)^2+2p^4(1-p)\left(1-p^2\right)+p^4\left(p(1-p)^2+(1-p)\right)\\
&+2p^5\left(1-p^2\right)+2p^6(1-p)+p^7 \\
=&p^3+3 p^4-2 p^5-2 p^6+p^7.
\end{aligned}
$$

For case (2), we have
$$
\begin{aligned}
\mathbb{E} (X_i X_j) & =2\left(p^3\left(p(1-p)^2+(1-p)\right)+p^4(1-p)^2+2p^5(1-p)\right)+p^6 \\
&=2 p^3+2 p^4-4 p^5+ p^6 .
\end{aligned}
$$

Note that $\mathbb{E}X_i\mathbb{E}X_j=(2p^2-p^4)^2$.
Let $f(p)=\mathbb{E}(X_iX_j)$
and $g(p)=\mathbb{E}X_i\mathbb{E}X_j+\delta(2p^2-p^4)^2$.
Then $f(1)=1<g(1)=1+\delta$.
Thus, for $p$ close enough to 1, we have $f(p)<g(p)$,
that is $\mathbb{E}( X_i X_j)-\mathbb{E}X_i\mathbb{E}X_j<\delta(2p^2-p^4)^2$.
Since there are at most $|\mathcal{F}_M|^2$ such pairs $\left(F_i, F_j\right)$,
the contribution of such pairs to $\operatorname{Var}_M X$ is at most
$|\mathcal{F}_M|^2\cdot \delta(2p^2-p^4)^2= \delta(\mathbb{E}X_M)^2.$

\vspace{2mm}
For (3), since all color classes in $\mathcal{F}_M$ are disjoint, for every $F_i \in\mathcal{F}_M$,
there are at most $4$ matchings $F_j$ that have at least an edge contained in $G [V(F_i)]$.
This means that there exist at most $2\cdot4 |\mathcal{F}_M|$ such pairs.
Since $\mathbb{E} (X_i X_j)-\mathbb{E} X_j \mathbb{E} X_j\leq 1$,
the contribution of such pairs to $\operatorname{Var} X_M$ is at most $O(n)$, which is $o((\mathbb{E}X_M)^2)$ as $\mathbb{E}X_M=\Omega(n)$.

\vspace{2mm}
Summing all above, we have
$$
0 \leq \operatorname{Var} X_M\leq \delta(\mathbb{E}X_M)^2+o\left((\mathbb{E} X_M)^2\right) .
$$
\end{proof}
By Inequality (\ref{e:1}),
for any $\delta>0$ we have
\begin{equation}\label{e:P1}
\mathbf{P}(\mathfrak{A}_M(\varepsilon))
\leq \mathbf{P}\left(X_{M} \leq (1-\varepsilon)\mathbb{E} X_{M}\right) \leq \frac{\operatorname{Var} X_{M}}{(\varepsilon \mathbb{E} X_{M})^2}=\dfrac{\delta}{{\varepsilon}^2}+o(1).
\end{equation}

\vspace{2mm}
\noindent{\bf Claim 2.} For any $\delta>0$,
we have
$0\leq \operatorname{Var} X_\Theta\leq \delta(\mathbb{E} X_{\Theta})^2$,
for $\Theta\in\{S,T\}$.

\begin{proof}
Take $\Theta=T$ as an example.
For any $F_i,F_j\in \mathcal{F}_T$,
if the edges in $F_i$ and $F_j$ are vertex-disjoint,
then $X_i$ and $X_j$ are independent and $\mathbb{E} (X_i X_j)-\mathbb{E} X_i \mathbb{E} X_j=0$.
Thus, assume that $F_i$ and $F_j$ are joint.
Then the contribution of $(F_i, F_j)$ to $\mathbb{E} (X_i X_j)$ is the probability that $X_i=X_j=1$.
Since $F_i$ and $F_j$ are edge-disjoint,
there is only one case if $F_i$ and $F_j$ are joint:
$F_i$ and $F_j$ intersect at exactly one common vertex.

Let $E(F_i)=\{a, b,c\}$
and  $E(F_j)=\{d,f,g\}$.
 Without loss of generality,
 we may assume that $V(a)\cap V(f)=\emptyset$.
Now we give the following table
in which $A$ denotes the event that the edges are in $G[H]$,
$B$ denotes the event that the edges are not in $G[H]$,
and
$P$ denotes the probability that
both the events $A$ and $B$ happen.

$$\begin{tabular}{|p{3cm}|p{4cm}|p{4cm}|}
\hline
$A$:&$B$:&$P(A\cap B)$\\
belong to $G[H]$&not belong to $G[H]$&\\
\hline
$a, f$&$b, c,d,h$&$p^4(1-p)$\\ \hline
$b, d $&$a, c,f,h$&$p^3(1-p)^2$\\ \hline
$b, h$&$a,c, d,f$&$p^3(1-p)^2$\\ \hline
$c, d $&$a,b ,f,h$&$p^3(1-p)^2$\\ \hline
$c,h$&$a,b,d,f$&$p^3(1-p)^2$\\ \hline
$a, b, c,d$&$f,h$&$p^4(1-p)$\\ \hline
$a,b,c,h$&$d,f$&$p^4(1-p)$\\ \hline
$d,h,f,b$&$a,c$&$p^4(1-p)$\\ \hline
$d,h,f,c$&$a,b$&$p^4(1-p)$\\ \hline
$a,b,c,d,h,f$&-&$p^5$\\ \hline
\end{tabular} $$

we omit the cases that do not exist in the table,  
such as
$a,c\in G[H]$ and $d,f,h,b\notin G[H]$.
Hence
$$\begin{aligned}
\mathbb{E}(X_iX_j)=5p^4(1-p)+4p^3(1-p)^2+p^5.
\end{aligned}$$

Note that $\mathbb{E}X_i\mathbb{E}X_j=(3p^2-2p^3)^2$.
Let $f(p)=\mathbb{E}(X_iX_j)$
and $g(p)=\mathbb{E}X_i\mathbb{E}X_j+\delta(3p^2-2p^3)^2$.
Then $f(1)=1<g(1)=1+\delta$.
Thus, for $p$ close enough to 1, we have $f(p)<g(p)$.
Since there are at most $|\mathcal{F}_T|^2$ such pairs $\left(F_i, F_j\right)$,
the contribution of such pairs to $\operatorname{Var} X_T$ is at most
$|\mathcal{F}_T|^2\cdot \delta(3p^2-2p^3)^2= \delta(\mathbb{E}X_T)^2.$

A similar argument as $\Theta=T$,
when $\Theta=S$
we also have
$$
0 \leq \operatorname{Var} X_S\leq \delta(\mathbb{E}X_S)^2 .
$$
\end{proof}

By Inequality (\ref{e:1}),
for any $\delta>0$ we have
\begin{equation}\label{e:P2}
\mathbf{P}(\mathfrak{A}_{\Theta}(\varepsilon))
\leq \mathbf{P}\left(X_{\Theta} \leq (1-\varepsilon)\mathbb{E} X_{\Theta}\right) \leq \frac{\operatorname{Var} X_{\Theta}}{(\varepsilon \mathbb{E} X_{\Theta})^2}=\dfrac{\delta}{{\varepsilon}^2},
\text{for~$\Theta\in \{T,S\}$}.
\end{equation}

For any $\varepsilon>0$,
let $\mathfrak{B(\varepsilon)}$ be the event that
$|H|\geq (1+\varepsilon)np$.
Since each vertex of $G$ is included in $H$ independently with probability $p$,
$|H|$ has the same probability distribution as $\operatorname{Bin}(n, p)$.
Then $\mathbb{E}|H|=np$.
Applying Chernoff's bound, we have
\begin{equation}\label{e:P3}
\mathbf{P}(\mathfrak{B(\varepsilon)})\leq \mathbf{P}\left(|H| \geq (1+\varepsilon)np\right) \leq \exp \left(-\frac{{\varepsilon}^2np}{3}\right).
\end{equation}

According to Inequalities (\ref{e:P1})-(\ref{e:P3}),
for any $\delta, \varepsilon>0$, we have
$$\begin{aligned}
\mathbf{P}(\mathfrak{A}_{M}(\varepsilon)
\cup\mathfrak{A}_{T}(\varepsilon)\cup\mathfrak{A}_{S}(\varepsilon)
\cup\mathfrak{B(\varepsilon)})&\leq \mathbf{P}(\mathfrak{A}_{M}(\varepsilon))+
\mathbf{P}(\mathfrak{A}_{T}(\varepsilon))+
\mathbf{P}(\mathfrak{A}_{S}(\varepsilon))+
\mathbf{P}(\mathfrak{B}(\varepsilon))\\
&\leq \frac{3\delta}{{\varepsilon}^2}+\exp \left(-\frac{{\varepsilon}^2np}{3}\right)+o(1).
\end{aligned}$$
Set $\delta\ll \varepsilon$,
there exists an integer $N_0$
such that when $n\geq N_0$,
$$\mathbf{P}(\mathfrak{A}_{M}(\varepsilon)
\cup\mathfrak{A}_{T}(\varepsilon)\cup\mathfrak{A}_{S}(\varepsilon)
\cup\mathfrak{B(\varepsilon)})<1,$$
which implies that there exists a subgraph $G[H]$ of $G$
satisfying:

(1) $|H| \leq (1+\varepsilon)n p$;

(2) the number of color classes
in $\mathcal{F}_{\Theta}$ that have at least one edge contained in $G[H]$ is at least $(1-\varepsilon)\cdot\mathbb{E}X_{\Theta}$
for $\Theta\in \{M,T,S\}$,
which implies that
the number of rainbow edges contained in $G[H]$ is at least
$
(1-\varepsilon)
\big(\left(2 p^2-p^4\right)\left|\mathcal{F}_M\right|+\left(3 p^2-2p^3\right)\left|\mathcal{F}_T\right| +p^2
\left|\mathcal{F}_S\right| \big) .
$

Let $f(p)=\frac{1}{n}\big(\left(2 p^2-p^4\right)\left|\mathcal{F}_M\right|+\left(3 p^2-2p^3\right)\left|\mathcal{F}_T\right| +p^2
\left|\mathcal{F}_S\right| \big)$ and $g(p)=p$.
Then $f'(p)=\frac{1}{n}\big((4p-4p^3)|\mathcal{F}_M|+(6p-6p^2)|\mathcal{F}_T|+2p|\mathcal{F}_S|\big)$
and $g'(p)=1$.
Since $\left|\mathcal{F}_M\right|+\left|\mathcal{F}_T\right|+\left|\mathcal{F}_S\right|=n$
and $\left|\mathcal{F}_M\right|+\left|\mathcal{F}_T\right|\geq \alpha n>\frac{n}{2}$,
we have $f(1)=g(1)=1$
and $f'(1)<1=g'(1)$.
Thus, for $p$ close enough to 1,
we have $f(p)>g(p)$.
Therefore, there exists a constant $\varepsilon>0$
such that
\begin{equation}\label{e:4}
(1-\varepsilon)
\big(\left(2 p^2-p^4\right)\left|\mathcal{F}_M\right|+\left(3 p^2-2p^3\right)\left|\mathcal{F}_T\right| +p^2
\left|\mathcal{F}_S\right| \big)>(1+\varepsilon)n p.
\end{equation}

Fixed $\varepsilon>0$ and $p$ closed enough to 1 satisfying Inequality (\ref{e:4}),
this completes the proof of
Theorem \ref{thm:main}  by setting $\beta=(1+\varepsilon)p$ and
$\gamma=\frac{1}{n}(1-\varepsilon)
\big(\left(2 p^2-p^4\right)\left|\mathcal{F}_M\right|+\left(3 p^2-2p^3\right)\left|\mathcal{F}_T\right| +p^2
\left|\mathcal{F}_S\right| \big)-(1+\varepsilon)p$
when $n\geq N_0$.
$\hfill\blacksquare$

\subsection{Proof of Theorem \ref{thm:main3}}

We prove Theorem \ref{thm:main3} by a random construction.
Let $G_0$ be an empty graph with $V(G_0)=\{v_1,v_2,\cdots,v_n\}$.
Let $\mathcal{V}_M$ be the set of all four-tuples in $V(G_0)$,
$\mathcal{V}_T$ be the set of all triples in $V(G_{0})$
and $\mathcal{V}_S$ be the set of all two-tuples of $V(G_0)$.
Let $\mathcal{V}=\mathcal{V}_T\cup \mathcal{V}_M \cup \mathcal{V}_S$.
Let $\mathcal{F}_{\Theta}$ be the set of elements in $\mathcal{V}_{\Theta}$
in which each element is included independently with probability $p=144n^{-3}$ for $\Theta\in\{M,T,S\}$. Let $\mathcal{F}=\mathcal{F}_M\cup \mathcal{F}_T \cup \mathcal{F}_S$.

\vspace{2mm}
\noindent{\bf Claim 1.} Let $\mathfrak{A}$ be the event that the number of pairs of distinct tuples in $\mathcal{F}$
that intersect with at least two vertices is at least $ n$.
Then $\mathbf{P}(\mathfrak{A})\leq \frac{6\cdot (144)^2}{ n}$.

\begin{proof}
For $V_i \in \mathcal{V}$,
let $X_i$ be the indicator random variable that $V_i \in \mathcal{F}$.
For the pair $(V_i, V_j) \in \left(\mathcal{F},\mathcal{F}\right)$,
let $X_{i,j}$ be the indicator random variable that $|V_i\cap V_j|\geq 2$.
It is clearly that $\mathbb{E} X_{i, j} \leqslant \mathbb{E}( X_i X_j)=p^2$ as $P(X_{i, j}=1) \leqslant P(X_i=X_j=1)=p^2$.
Let $X=\sum_{(V_i, V_j) \in \left(\mathcal{F},\mathcal{F}\right)}X_{i,j}$.

Note that there are
at most $\binom {n}{3} \cdot 3 \cdot(n-2)$ pairs $(V_i, V_j) \in(\mathcal{V}_T, \mathcal{V}_T)$
such that $|V_i\cap V_j|\geq 2$
(as there are at most $\binom {n}{3}$ ways to choose $V_i \in \mathcal{V}_T$,
3 ways to choose two intersecting vertices,
and then at most $(n-2)$ ways to choose the last vertex of $V_j$)
and
at most $\binom{n}{4}\cdot\binom{4}{2} \cdot(n-2)$ pairs $(V_i, V_j) \in(\mathcal{V}_T, \mathcal{V}_M)$
such that $|V_i\cap V_j|\geq 2$
(as there are at most $\binom{n}{4}$ ways to choose $V_i \in \mathcal{V}_M$,
$\binom{4}{2}$ ways to choose two intersecting vertices,
and then at most $(n-2)$ ways to choose the last vertex of $V_j$).
Similarly, there are at most $\binom {n}{3} \cdot 3$  pairs $(V_i, V_j) \in(\mathcal{V}_T, \mathcal{V}_S)$,
$\binom{n}{4}\cdot\binom{4}{2} \cdot\binom{n-2}{2}$ pairs $(V_i, V_j) \in(\mathcal{V}_M, \mathcal{V}_M)$,
$\binom{n}{4}\cdot\binom{4}{2}$ pairs $(V_i, V_j) \in(\mathcal{V}_M, \mathcal{V}_S)$,
and $\binom{n}{2}$ pairs $(V_i, V_j) \in(\mathcal{V}_S, \mathcal{V}_S)$
satisfying $|V_i\cap V_j|\geq 2$.
Therefore,

$$\begin{aligned}
\mathbb{E}X=&\mathbb{E}(\sum_{F_i, F_j \in \left(\mathcal{F},\mathcal{F}\right)}X_{i,j})\\
&\leq \big(\binom {n}{3} \cdot 3 \cdot(n-2)+\binom{n}{4}\cdot\binom{4}{2} \cdot(n-2)+\binom {n}{3} \cdot 3\\
&+\binom{n}{4}\cdot\binom{4}{2} \cdot\binom{n-2}{2}+\binom{n}{4}\cdot\binom{4}{2}+\binom{n}{2}\big)\cdot p^2\\
\leq& 6n^6p^2.
\end{aligned}$$

By Markov's inequality, we have $\mathbf{P}(X\geq n)\leq \frac{6\cdot (144)^2}{ n}$.
\end{proof}

Now we define a new system $E(\mathcal{F})$ consisting of  $e\in \mathcal{V}_S$
such that for each $e\in E(\mathcal{F})$
there exists an element $V\in \mathcal{F}$ with
$e\subseteq V$.
A cycle $C=(e_1,e_2,\cdots,e_\ell)$ of length $\ell$ in $E(\mathcal{F})$ is a family of $\ell$ elements in $E(\mathcal{F})$
in which $e_i$ is intersect with successor $e_{i+1}$ and predecessor $e_{i-1}$ on $C$ at distinct vertices
for $i\in [\ell]$.
We call $C$ is \emph{rainbow} if
for $i\in [\ell]$, there exists an element $V\in \mathcal{F}$ 
such that $e_i\subseteq V$ and $e_j\nsubseteq V$
for each $j\in  [\ell]\setminus \{i\}$.

\vspace{2mm}
\noindent{\bf Claim 2.} Let $\mathfrak{B}$ be the event that the number of rainbow cycles of length
at most $c\log n$ is at least $ n$.
Then $\mathbf{P}(\mathfrak{B})\leq n^{-\frac{1}{2}}\cdot c\log n$.

\begin{proof}
For a fixed 2-tuple $e\in \mathcal{V}_S$,
the probability that there exists no $V\in \mathcal{F}$ such that $e\subseteq V$
is at least $(1-p)^{1+n-2+\binom{n-2}{2}}$
as we can find one tuple in $\mathcal{V}_S$, $n-2$ tuples in $\mathcal{V}_T$,
and $\binom{n-2}{2}$ tuples in $\mathcal{V}_M$ that contain $e$.
Hence we have  $$\begin{aligned}
\mathbf{P}(e\in E(\mathcal{F}))&=\mathbf{P}(e\subseteq V \text { for some } V \in \mathcal{F})\\
&=
1-\mathbf{P}(e\nsubseteq F \text { for any } F \in \mathcal{F})\\
&\leq1-(1-p)^{(1+n-2+\binom{n-2}{2})}\\
&\leq \left(1+n-2+\binom{n-2}{2}\right)p.
\end{aligned}
$$

Let $C_k=(e_1,e_2,\cdots,e_k)$ be a cycle of length $k$ in $\mathcal{F}$.
For $1\leq i\leq k$, we define $\mathfrak{S}_{C}(e_i)$ be the event that
there exists an $V\in \mathcal{F}$ such that
$e_i\subseteq V$ while $e_j\nsubseteq V$ for $j<i$ on $C$.
Since there are at most $1+n-2+\binom{n-2}{2}$
tuples in $\mathcal{V}$ containing $e_i$,
we have $\mathbf{P}\left(\overline{\mathfrak{S}_{C}(e_i)} \mid \cap_{j=1}^{i-1} \mathfrak{S}_{C}(e_j)\right) \geq(1-p)^{1+n-2+\binom{n-2}{2}}$,
as
we can assume that these tuples are not in $\mathcal{F}$.
Therefore,
$$
\begin{aligned}
\mathbf{P}\left(C_k  \text { is rainbow }\right)
&\leq \mathbf{P}\left(\bigcap_{i=1}^k \mathfrak{S}_{C}(e_i)\right)=\prod_{i=1}^k \mathbf{P}\left(\mathfrak{S}_{C}(e_i) \mid \bigcap_{j=1}^{i-1} \mathfrak{S}_{C}(e_j)\right)\\
&\leq \big(1-(1-p)^{1+n-2+\binom{n-2}{2}}\big)^k\\
&\leq\left(1+n-2+\binom{n-2}{2}\right)^kp^k.
\end{aligned}
$$

Let $Y_k$ be the number of rainbow cycles of length $k$ in $E(\mathcal{F})$.
Since any consecutive tuples in $C_k$ are intersect,
the number of cycles of length $k$ in $E(\mathcal{F})$ is at most $n^k$.
Then
$$
\mathbb{E} Y_k \leq n^k \left(1+n-2+\binom{n-2}{2}\right)^kp^k \leq n^{1 / 2}
$$
for $k \leq c\log n$ and $c>0$ small enough.
Let $Y=\sum_{k=3}^{\lfloor c\log n \rfloor}Y_k$.
Therefore
$$
\mathbb{E} Y=\sum_{k=3}^{\lfloor c \log n\rfloor} \mathbb{E} Y_k\leq n^{\frac{1}{2}}\cdot c\log n .
$$
Thus by Markov's inequality, we have
$$
\mathbf{P}\left(Y \geq n\right)\leq n^{-\frac{1}{2}}\cdot c\log n.
$$
\end{proof}

\noindent{\bf Claim 3.} Let $\mathfrak{C}$ be the event that $|\mathcal{F}|\leq 3 n$.
Then $\mathbf{P}(\mathfrak{C})\leq \exp(-\frac{4n}{200})$.

\begin{proof}
Since each tuple in $\mathcal{V}$ is included in $\mathcal{F}$ independently with probability $p=144n^{-3}$,
we have
$$
\mathbb{E}|\mathcal{F}|
=\mathbb{E}|\mathcal{F}_{M}|+\mathbb{E}|\mathcal{F}_{T}|+\mathbb{E}|\mathcal{F}_{S}|
=\left(\binom{n}{2}+\binom{n}{3}+\binom{n}{4}\right) p \geq 4  n .
$$
Applying Chernoff's bound, we have
$$
\mathbf{P}(|\mathcal{F}| \leq 3  n) \leq \mathbf{P}(|\mathcal{F}| \leq 0.9 \cdot \mathbb{E}|\mathcal{F}|)\leq \exp(-\frac{4n}{200}).
$$
\end{proof}

According to Claims 1-3,
we have
$$\begin{aligned}
\mathbf{P}(\mathfrak{A}
\cup\mathfrak{B}\cup\mathfrak{C})
&\leq \mathbf{P}(\mathfrak{A})+\mathbf{P}(\mathfrak{B})+\mathbf{P}(\mathfrak{C})\\
&\leq \frac{6}{n}+n^{-\frac{1}{2}}\cdot c\log n+\exp \left(-\frac{4n}{200}\right).
\end{aligned}$$
Thus
there exists an integer $N_0$
such that when $n\geq N_0$,
$$\mathbf{P}(\mathfrak{A}
\cup\mathfrak{B}\cup\mathfrak{C})<1,$$
which implies that there exists a  system $\mathcal{F}$ including
2-tuples, triples and 4-tuples satisfying

(1) $|\mathcal{F}|\geq 3 n$;

(2) the number of pairs of distinct tuples in $\mathcal{F}$
that intersect with at least two vertices is at most $ n$;

(3) the number of rainbow cycles of length
at most $c\log n$ is at most $n$.

Now, we can remove one tuple in the pairs that intersect with at least two vertices from $\mathcal{F}$ to get $\mathcal{F}'$.
Then the tuples in $\mathcal{F}'$ intersect with each other at most one vertex.
In particular, each $e \in E\left(\mathcal{F}\right)$ is contained in exactly one $F \in \mathcal{F}'$. Then
$$
\left|\mathcal{F}'\right| \geq 3 n- n \geq 2  n \text {. }
$$

Next, for each rainbow cycle of length at most $c \log n$ in $\mathcal{F}'$,
we can choose one edge $e$ and remove the tuple in $\mathcal{F}'$ that contains $e$ to get $\mathcal{F}''$.
In particular, each rainbow cycle in $\mathcal{F}''$ is
of length more than $c \log n$.
Therefore
$$
\left|\mathcal{F}''\right| \geq 2  n- n \geq  n \text {. }
$$

Let $G:=(V(G_0),E\left(\mathcal{F}''\right))$. Then $G$ is a graph formed by at least $n$ edge-disjoint color classes,
in which each color class is a
matching or triangle or single edge,
and without rainbow cycles of length less than $c \log n$. This completes the proof.
$\hfill\blacksquare$

\section*{Acknowledgment}
The first author thanks Bo Ning for giving the problem studied in this paper to her.
This research was supported in part by the National Natural Science Foundation of China under grant numbers 12301457 and 12171440.

\end{document}